\documentclass{article}[12pt]

\usepackage{arxiv}

\usepackage[utf8]{inputenc} % allow utf-8 input
\usepackage[T1]{fontenc}    % use 8-bit T1 fonts
\usepackage{hyperref}       % hyperlinks
\usepackage{url}            % simple URL typesetting
\usepackage{booktabs}       % professional-quality tables
\usepackage{amsfonts}       % blackboard math symbols
\usepackage{nicefrac}       % compact symbols for 1/2, etc.
\usepackage{microtype}      % microtypography
\usepackage{lipsum}
\usepackage{graphicx}
\title{On Dropping Needles and WiFi Link Crossing}

\author{
  Anurag Pallaprolu\thanks{https://www.github.com/apallaprolu} \\
  Electronics and Computer Engineering\\
  University of California, Santa Barbara\\
  Santa Barbara, CA 93106 \\
  \texttt{apallaprolu@ucsb.edu} \\
}

\begin{document}
\maketitle

\begin{abstract}
In a general simulation of random walking (with the angle of motion picked uniformly), it can be seen that the probability of crossing a WiFi TX-RX link is directly proportional to the per-step distance and inversely proportional to the lateral dimension of the room. The asymptotic value of the said crossing probability is derived using Perron-Frobenius theory to determine the limit distribution of the said Markov model. Surprisingly, we can establish a bijection to a scenario explored nearly 300 years ago by Georges-Louis Leclerc, Comte de Buffon to get the result. Furthermore we can use the generalizations of the latter problem to ascertain some interesting observations about the original one. 
\end{abstract}

% keywords can be removed

\section{Introduction}
It has been definitively shown in \cite{mostofi1} and \cite{mostofi2} that for Markov walking in a rectilinear enclosure of length $L$ the limit distributions for the position and angular orientation are uniform with the probability of crossing a link placed at $(\frac{L}{2}, 0), (\frac{L}{2}, B)$ given by
$$p_c = \frac{2d_s}{\pi L}$$
For the sake of completeness, we will state below the iteration steps for simulating this model:
\begin{itemize}
    \item Set initial angle $\theta_0$ and initial coordinates $x_0, y_0$.
    \item Update the next angle by picking $\theta_{i}$:
    \begin{itemize}
        \item $\theta_{i} = \theta_{i-1}$ with probability $p_\theta$
        \item $\theta_{i} \in [\alpha_i]$ with probability $1 - p_\theta$ where the $\alpha$ set is the uniform distribution over a certain allowed angular range.
    \end{itemize}
    \item $x_i = x_{i-1} + d_s\cos\theta_i$
    \item $y_i = y_{i-1} + d_s\sin\theta_i$
    \item Collisions with the boundary are handled by assuming that there is a perfect reflection upon contact. 
\end{itemize}
In this paper we will try to connect this model with Buffon's needle problem \cite{buffon1}, whose statement below is adapted from \cite{krantz}):
\begin{quote}
    A floor is made of long planks that are $L$ inches wide. A girl drops a thin stick that is $d_s$ inches long onto the floor. She does so a great many times- say N times. Calculate the probability that the stick will land on a crack between the planks in the floor as a function of $N$. What is the asymptotic value as $N\to\infty$?
\end{quote}
\section{The Connection}
The solution is rather straightforward if we first calculate the joint density function of the distance and angle of landing. If we assume that the midpoint of stick is $d_m$ units away from a crack and if the extrapolated angle between the stick and the crack is $\theta$, we have a crossing if
$$0 \leq d_m \leq \frac{d_s}{2}\sin\theta$$
\begin{figure}
    \centering
    \includegraphics[scale=0.6]{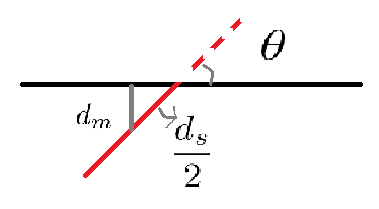}
    \includegraphics[scale=0.6]{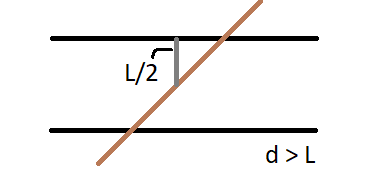}
\end{figure}
The probability density for $d_m$ consequently is uniform with value $\frac{2}{L}$ since the maximum length of a stick that has a non-zero chance of not crossing is nothing but the spacing between the lines. Also, the angular variable has a density of $\frac{2}{\pi}$ since the stick can go all the way from being parallel to perpendicular to the cracks. We can therefore write the joint probability density as
\\
$$f(d_m, \theta) = \frac{4}{\pi L}$$
\\
To find the probability of crossing, we just integrate the density over the reduced sample space represented by the bound established before:
\\
$$P = \int_{0}^{\frac{\pi}{2}}\int_{0}^{\frac{d_s}{2}\sin\theta} \frac{4}{\pi L} d{d_m}d\theta = \frac{2d_s}{\pi L}$$
\\
From the figure (which shows only one crack), we can infer that this analysis is valid for $d_s << L$ since if it is not the case, we would have $d_m = \frac{L}{2}$ for all angles where there is a crossing on two consecutive cracks, and the analysis gets further involving.
\\
\\
To show the relation between this analysis to that of the probability of crossing a WiFi link, we construct an ensemble of rectilinear rooms of length L, each having a TX-RX link at $\frac{L}{2}$ abutted length-wise, as shown below:
\begin{center}
\includegraphics[scale=0.6]{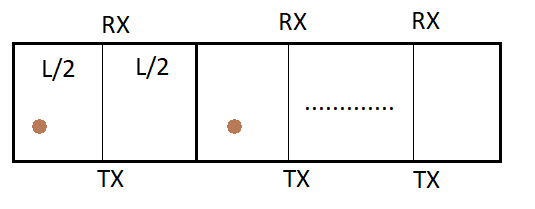}
\end{center}    
It should be rather obvious that a link crossing would take place if statistically we can see that a metaphorical "stick" with length equal to the size of a single step of the walk happens to fall on the link which resembles the "crack" in the earlier discussion. Since our model assumes the reflection policy for boundary collisions, we have the symmetric behavior for crossing approached from $0<x_l<\frac{L}{2}$ or $\frac{L}{2}<x_l<L$ where $x_l$ is the location of the link. Indeed, this is the same result that is obtained after the much more mathematically involving analysis presented in \cite{mostofi1}.
\newpage
\section{Buffon's Noodle}
Boris Gnedenko proves in his classic text \cite{gnedenko} that our discussion above need not be restricted to movements along a straight line. In-fact, if we have any \textit{rectifiable} curve of motion, the probability of crossing remains invariant. This idea is colloquially known as Buffon's Noodle problem instead, although this result generalizes the former. For the sake of completeness, we will state the question and a sketch of the proof for the same below.
\begin{quote}
    Throw at random a convex contour of diameter less than $L$ on a horizontal plane partitioned by parallel lines separated by a distance $L$. Find the probability that the contour will intersect a line. 
\end{quote}
The solution would involve splitting the curve into segments and applying the result of the earlier section for each of them. The requirement of the curve being rectifiable pops up here since we would need the line integral along the curve to be bounded (which fails for some pathological examples). Let us assume the split
$$\mathcal{P} = [p_{12}, p_{23}, p_{34}, ....p_{n-1, n}]$$
We must emphasize that this split/rectification is a special one where each $p_{ij} << L$ so the cardinality of $\mathcal{P}$ would be large. We can apply our earlier result to see that $p_c = \frac{2p_{ij}}{\pi L}$ for each segment. If we assume each step in the path is of independent movements (the validity of the reduction of every action in the sample space to a basis of independent events is a controversial point in the theory of geometric probability as explained in \cite{gnedenko}) then we simply sum the probabilities of crossing for each segment and we obtain the desired result. However, if we are interested in find the average number of crossings when a convex contour of length $d_s$ is thrown on the parallel line lattice, we do not need the independence condition. By directly casting the law of linearity of expectation (inspired from \cite{ramaley}) we can see that the expected number of crossings will be invariant in general.  

\section{Conclusions}
In this document, we have shown a much simpler proof of the probability of crossing a WiFi link in a room as a pure function of the length measures. In doing so we have shown that the result holds good even if the location of the link is not fixed and also if the movement is not precisely rectilinear. 

\bibliographystyle{unsrt}  


\end{document}